\documentclass[12pt,A4]{article}
\usepackage[utf8]{inputenc}
\usepackage{amsmath}
\usepackage{commath}
\usepackage{authblk}

\title{Exponentiation Using Laplace Expansion}

\author{Bhavesh Lakhotia}

\affil{Institute and Faculty of Actuaries, UK}
\date{15 June 2021}

\begin{document}

\maketitle

\section{Introduction}

The Laplace expansion is a well known method for calculating determinants of matrices \cite{Laplace}. In this paper, I use the Laplace expansion to derive an equation. The equation essentially allows for the exponentiation of a base, given an integer exponent. It can be used to calculate exponents using parallel processors. To the best of my knowledge, this method of exponentiation is not being currently used in the literature.

\section{Define P}
Let $a \neq P_1 \neq P_2 \neq P_3 \neq P_4 \cdots P_k$ be $k+1$ elements of a field $F$. Define the $k \times k$ matrix P = \Large{\[
 \begin{bmatrix}
  P^{0}_{1} & P^{1}_{1} & P^{2}_{1}  &\cdots& P^{k-2}_{1}& \sum^{k-1}_{i=1} P^{(k-1)-i}_{1}a^{i-1}\\
  P^{0}_{2} & P^{1}_{2} & P^{2}_{2}  &\cdots& P^{k-2}_{2}&\sum^{k-1}_{i=1} P^{(k-1)-i}_{2}a^{i-1}\\
  P^{0}_{3} & P^{1}_{3} & P^{2}_{3}  &\cdots& P^{k-2}_{3}&\sum^{k-1}_{i=1} P^{(k-1)-i}_{3}a^{i-1}\\
  \vdots & \vdots & \vdots & \cdots & \vdots &\vdots \\
  P^{0}_{k} & P^{1}_{k} & P^{2}_{k} & \cdots & P^{k-2}_{k}& \sum^{k-1}_{i=1} P^{(k-1)-i}_{k}a^{i-1}\\
 \end{bmatrix}
\]}
\\

This matrix has zero determinant as column $k$ is a linear combinations of the other columns. We can rewrite the matrix as below.\\ 

\Large{\[
 \begin{bmatrix}
  P^{0}_{1} & P^{1}_{1} & P^{2}_{1}  &\cdots& P^{k-2}_{1}& \frac{P^{k-1}_{1} - a^{k-1}}{P_{1}- a}\\
  P^{0}_{2} & P^{1}_{2} & P^{2}_{2}  &\cdots& P^{k-2}_{2}& \frac{P^{k-1}_{2} - a^{k-1}}{P_{2}- a}\\
  P^{0}_{3} & P^{1}_{3} & P^{2}_{3}  &\cdots& P^{k-2}_{3}&\frac{P^{k-1}_{3} - a^{k-1}}{P_{3}- a}\\
  \vdots & \vdots & \vdots & \cdots & \vdots &\vdots \\
  P^{0}_{k} & P^{1}_{k} & P^{2}_{k} & \cdots & P^{k-2}_{k}& \frac{P^{k-1}_{k} - a^{k-1}}{P_{k}- a}\\
 \end{bmatrix}
\]}

\section{Define V}
Let V(P) be the Vandermonde matrix associated P. V(P) is thus defined:\\ 

\Large{\[
 \begin{bmatrix}
  P^{0}_{1} & P^{1}_{1} & P^{2}_{1}  &\cdots& P^{k-2}_{1}& P^{k-1}_{1} \\
  P^{0}_{2} & P^{1}_{2} & P^{2}_{2}  &\cdots& P^{k-2}_{2}& P^{k-1}_{2} \\
  P^{0}_{3} & P^{1}_{3} & P^{2}_{3}  &\cdots& P^{k-2}_{3}&P^{k-1}_{3} \\
  \vdots & \vdots & \vdots & \cdots & \vdots &\vdots \\
  P^{0}_{k} & P^{1}_{k} & P^{2}_{k} & \cdots & P^{k-2}_{k}& P^{k-1}_{k}\\
 \end{bmatrix}
\]}

The cofactors of the $k^{th}$ column are the same for V(P) and P. The determinant of V(P) is non-zero and given by \norm{V(P)}. It is known that \cite{Turner}:
\\
\begin{align}
\norm{V(P)} = \LARGE{\prod_{1\leq i<j\leq k} (P_{j} -P_{i}) }
\end{align}

\section{Characterizing P}

Let $C_{j,k}$ denote the $j^{th}$ cofactor of the $k_{th}$ column of P. Then, using the Laplace expansion of the $k^{th}$ column we have,

\begin{align}
\LARGE{\sum_{j=1}^{k}\frac{(P^{k-1}_{j} - a^{k-1}) \times C_{j,k}}{P_{j} - a} = 0 }
\end{align}

Rearranging the terms, we have,
\\

\begin{align}
\LARGE{\sum_{j=1}^{k} \frac{P^{k-1}_{j}  \times C_{j,k}}{P_{j} - a} = (a^{k-1}) \times (\sum_{j=1}^{k} \frac{ C_{j,k}}{P_{j} - a}) }
\end{align}

In the equation above, $(\sum_{j=1}^{k}\frac{ C_{j,k}}{P_{j} - a})$ is non-zero. See appendix for this. Thus, we have the final equation for computing the $(k-1)^{th}$ power of a: \\
\begin{equation}
 \left[  \frac{  \sum_{j=1}^{k} \frac{P^{k-1}_{j} \times C_{j,k}}{P_{j} - a} } { \sum_{j=1}^{k} \frac{ C_{j,k}}{P_{j} - a} }  \right] = a^{k-1} \ \ \ \ 
 \end{equation}\\

Note that in equation (1), $C_{j,k}$  has the form\\
$\frac{ (-1)^{k+j} \times  \norm{V(P)}}{\prod_{1\leq i<j} (P_{j} -P_{i}) \times \prod_{j< i \leq k } (P_{i} -P_{j})} $\\

This is because the the minor $M_{j,k}$ is of the form:\\
$\frac{  \norm{V(P)}}{\prod_{1\leq i<j} (P_{j} -P_{i}) \times \prod_{j< i \leq k } (P_{i} -P_{j})} $\\

Define $C_{j} = \frac{1}{ \prod_{1 \leq i \leq k, j \neq i }(P_{i} -P_{j})} $. Then we have the simplified equation,\\
\begin{equation}
 \left[  \frac{  \sum_{j=1}^{k} \frac{P^{k-1}_{j} \times C_{j}}{P_{j} - a} } { \sum_{j=1}^{k} \frac{ C_{j}}{P_{j} - a} }  \right] = a^{k-1} \ \ \ \ 
 \end{equation}\\

This is because,\\
\\
$ (-1)^{k+j} \times  \prod_{1\leq i<j} (P_{j} -P_{i}) =  (-1)^{k} \times \prod_{1\leq i<j} (P_{i} -P_{j}) $\\

$(-1)^k$ and $\norm{V}$ are common to both the numerator and denominator in equation (1) so can be ignored when defining $C_j$.\\

Equation (2) is the central result of this paper. In the equation, the $C_j$ are scale and shift invariant with respect to $P_{j}$ after cancellations. This means that for a given $ \alpha$ and $\beta \neq 0$ in $F$,  $ a \neq Q_j = \alpha + \beta \times P_{j}$ $ \forall$ $ 1 \leq j \leq k $ satisfies the equation,\\

\begin{equation}
 \left[  \frac{  \sum_{j=1}^{k} \frac{Q^{k-1}_{j} \times C_{j}}{Q_{j} - a} } { \sum_{j=1}^{k} \frac{ C_{j}}{Q_{j} - a} }  \right] = a^{k-1} \ \ \ \ 
 \end{equation}\\

It can be shown that the following identity, which is an extension of (2), also holds for an integer $n$, $0 \leq n<k$. To see this, simply replace the $k^{th}$ column of P with the ${n+1}^{th}$ column of $P$ for $0 \leq n< (k-1)$. Then take the Laplace expansion along column $k$.  Recognize this determinant to be 0 as two columns of the resulting matrix are equal. Thus,  $\sum_{1 \leq j \leq k} P^n_j \times C_{j,k} =  0 = \sum_{1 \leq j \leq k} P^n_j \times C_{j}  $.\\

\begin{equation}
 \left[  \frac{  \sum_{j=1}^{k} \frac{P^{n}_{j} \times C_{j}}{P_{j} - a} } { \sum_{j=1}^{k} \frac{ C_{j}}{P_{j} - a} }  \right] = a^{n} \ \ \ \ 
 \end{equation}\\

 \section{Discussion}

In this section, I discuss the possibility of faster exponentiation using  equation (2). In equation(2), depending on the choice of $P_{j}$, the equation can take interesting forms. Two forms are discussed below. In all cases the non-equality assumption holds for the elements $P_j$ and base $a$. Note also that $P^{k-1}_j$ and  $C_j$ do not depend on $a$. They can be pre-calculated and stored in processors. Suppose these processors are connected in parallel.  Then, first, each summand in the numerator and denominator in equation (2) can be calculated simultaneously once $a$ is known. Second, the summation to get the numerator and to get the denominator can then be made in close to $ \log_{2} k$ steps. Third, a final division gives the requisite exponentiation. The speed of these 3 steps determines the speed of computation. Note that the possibility of faster exponentiation arises because addition is typically faster than multiplication in the second step.

As a first example, let $C_j = j$ in equation (2). Then, multiplying the numerator and denominator by $(k-1)!$ gives us this form of the equation.\\
 \begin{equation}
 \left[  \frac{  \sum_{j=1}^{k} \frac{j^{k-1} \times (-1)^{j} \times \binom{k-1}{j-1}}{j - a} } { \sum_{j=1}^{k} \frac{ (-1)^{j} \times \binom{k-1}{j-1}}{j - a} }  \right] = (a)^{k-1} \ \ \ \ 
  \end{equation}

 Second, let  $P_j = e_j$ in equation (2), where $e_j$ are the $(k-1)^{th}$ roots of unity such that $e_1 \neq e_2 \neq e_3 \neq \cdots \neq e_{k-1} $. Then, choose $P_k$ of a choice that assists computation. For example, choose $P_k =0$. Now, subtract 1 from both sides of the equation to get:\\

 \begin{equation}
 \left[  \frac{   \frac{( - 1) \times C_k  }{- a} } { (\sum_{j=1}^{k-1} \frac{  C_j}{e_j - a} ) + \frac{C_k}{-a}}  \right] = (a)^{k-1}-1\ \ \ \ 
  \end{equation}\\
 This form can further be simplified to, 
 
 \begin{equation}
 \left[  \frac{   ( - 1) } { \sum_{j=1}^{k-1} \frac{ (- a) \times C_j}{(e_j - a) \times C_k} +1 }  \right] = (a)^{k-1}-1\ \ \ \ 
  \end{equation}

 Recognising further that $\sum_{1\leq i \leq k} C_i = 0$, this is further simplified to,

 \begin{equation}
 \left[  \frac{   -1 } { \sum_{j=1}^{k-1} \frac{ (- e_j) \times C_j}{(e_j - a) \times C_k}  }  \right] = (a)^{k-1}-1\ \ 
  \end{equation}
 
 $C_k=(-1)^{k}$ since it is the product of all distinct $(k-1)^{th}$ roots of unity, we have,\\ 
 
  \begin{equation}
  \left[  \frac{   (-1)^{k} } { \sum_{j=1}^{k-1} \frac{ (e_j) \times C_j}{(e_j - a)}  }  \right] = (a)^{k-1}-1\ \ 
  \end{equation}
 
 When this equation is being used for exponentiation in a computer, this essentially requires the calculation of $\frac{1}{e_j -a}$ in each of the (k-1) nodes of the parallel processor. This is because all other values on the LHS of the equation are independent of $a$ and can be pre-calculated. 
 
 Equation(9) is closely related to the following identity \cite{RootsUnity}.\\
 \begin{equation}
 {(a)^{k-1}}-1 = \LARGE{\prod_{1\leq i\leq (k-1)} (a -e_{i}) }
 \end{equation}\\
 
Equation(10) can be rewritten as:\\
 \begin{equation}
 \begin{split}
 \frac{1}{(a)^{k-1}-1} & = \frac{1}{\LARGE{\prod_{1\leq i\leq (k-1)} (a -e_{i}) }} \\
 &= \sum_{1 \leq i \leq (k-1)} \frac{1}{\LARGE{(a-e_i) \times \prod_{1\leq j\leq (k-1), i \neq j} (e_i -e_{j}) }}
 \end{split}
  \end{equation}

 The equivalence of equation (9) and equation (11) is not difficult to see.\\

As shown in this paper, interesting forms can be derived with the use of equation(2). The use of equation (2) using a parallel computing structure is the subject of the patent filed with the patent office in India with the title "Exponential Calculator Using Parallel Processor Systems" with application number 202111023892. This patent discloses how exponentiation can be done efficiently using the ideas of this paper.

\section{Appendix}

Let P and V(P) be defined as in section 2 and section 3 respectively. Also, let $c = \prod_{1\leq j\leq k} (P_{j} -a)^{-1} \neq 0$. The purpose of this appendix is to show $(\sum_{j=1}^{k}\frac{ C_{j,k}}{P_{j} - a})$ is non-zero. To see this, note that this expression is the determinant of X, defined below.

Define  $X= $\\
\Large{\[
 \begin{bmatrix}
  P^{0}_{1} & P^{1}_{1} & P^{2}_{1}  &\cdots& P^{k-2}_{1}& \frac{1}{P_{1} -a} \\
  P^{0}_{2} & P^{1}_{2} & P^{2}_{2}  &\cdots& P^{k-2}_{2}& \frac{1}{P_{2} -a} \\
  P^{0}_{3} & P^{1}_{3} & P^{2}_{3}  &\cdots& P^{k-2}_{3}& \frac{1}{P_{3} -a} \\
  \vdots & \vdots & \vdots & \cdots & \vdots &\vdots \\
  P^{0}_{k} & P^{1}_{k} & P^{2}_{k} & \cdots & P^{k-2}_{k}& \frac{1}{P_{k} -a}\\
 \end{bmatrix}
\]}

Then we have,
\begin{equation*}
\norm{X}= \Large{
 \begin{vmatrix}
  P^{0}_{1} & P^{1}_{1} & P^{2}_{1}  &\cdots& P^{k-2}_{1}& \frac{1}{P_{1} -a} \\
  P^{0}_{2} & P^{1}_{2} & P^{2}_{2}  &\cdots& P^{k-2}_{2}& \frac{1}{P_{2} -a} \\
  P^{0}_{3} & P^{1}_{3} & P^{2}_{3}  &\cdots& P^{k-2}_{3}& \frac{1}{P_{3} -a} \\
  \vdots & \vdots & \vdots & \cdots & \vdots &\vdots \\
  P^{0}_{k} & P^{1}_{k} & P^{2}_{k} & \cdots & P^{k-2}_{k}& \frac{1}{P_{k} -a}\\
 \end{vmatrix}
}
  \end{equation*}

Since $m^{0} = 1$,
\begin{equation*}
\norm{X}=\Large{
 \begin{vmatrix}
  1 & P^{1}_{1} & P^{2}_{1}  &\cdots& P^{k-2}_{1}& \frac{1}{P_{1} -a} \\
  1 & P^{1}_{2} & P^{2}_{2}  &\cdots& P^{k-2}_{2}& \frac{1}{P_{2} -a} \\
  1 & P^{1}_{3} & P^{2}_{3}  &\cdots& P^{k-2}_{3}& \frac{1}{P_{3} -a} \\
  \vdots & \vdots & \vdots & \cdots & \vdots &\vdots \\
  1 & P^{1}_{k} & P^{2}_{k} & \cdots & P^{k-2}_{k}& \frac{1}{P_{k} -a}\\
 \end{vmatrix}
}
  \end{equation*}

Multiply RHS by $c$. Multiply $j^{th}$ row by $P_{j}-a$ to get,
\begin{equation*}
c \times
\Large{
 \begin{vmatrix}
  P^{1}_{1} - a  & P^{2}_{1} - a P^{1}_{1} & P^{3}_{1} - a P^{2}_{1}  &\cdots& P^{k-1}_{1} - a P^{k-2}_{1}& 1\\
  P^{1}_{2} - a  & P^{2}_{2} - a P^{1}_{2} & P^{3}_{2} - a P^{2}_{2}  &\cdots& P^{k-1}_{2} - a P^{k-2}_{2}& 1\\
  P^{1}_{3} - a  & P^{2}_{3} - a P^{1}_{3} & P^{3}_{3} - a P^{2}_{3}  &\cdots& P^{k-1}_{3} - a P^{k-2}_{3}& 1\\
  \vdots & \vdots & \vdots & \cdots & \vdots &\vdots \\
  P^{1}_{k} - a  & P^{2}_{k} - a P^{1}_{k} & P^{3}_{k} - a P^{2}_{k}  &\cdots& P^{k-1}_{k} - a P^{k-2}_{k}& 1 \\
 \end{vmatrix}
}\\
 \end{equation*}
 
Multiply column 1 by $a$ and add to column 2 to get,
\begin{equation*}
c \times
\Large{
 \begin{vmatrix}
  P^{1}_{1} - a  & P^{2}_{1} - a^{2}  & P^{3}_{1} - a P^{2}_{1}  &\cdots& P^{k-1}_{1} - a P^{k-2}_{1}& 1 \\
  P^{1}_{2} - a  & P^{2}_{2} - a^{2}  & P^{3}_{2} - a P^{2}_{2}  &\cdots& P^{k-1}_{2} - a P^{k-2}_{2}& 1 \\
  P^{1}_{3} - a  & P^{2}_{3} - a^{2}  & P^{3}_{3} - a P^{2}_{3}  &\cdots& P^{k-1}_{3} - a P^{k-2}_{3}& 1 \\
  \vdots & \vdots & \vdots & \cdots & \vdots &\vdots \\
  P^{1}_{k} - a  & P^{2}_{k} - a^{2}  & P^{3}_{k} - a P^{2}_{k}  &\cdots& P^{k-1}_{k} - a P^{k-2}_{k}& 1 \\
 \end{vmatrix}
}\\
\end{equation*}

Multiply column 2 by $a$ and add to column 3 to get,
\begin{equation*}
c \times
\Large{
 \begin{vmatrix}
  P^{1}_{1} - a  & P^{2}_{1} - a^{2}  & P^{3}_{1} - a^{3}   &\cdots& P^{k-1}_{1} - a P^{k-2}_{1}& 1 \\
  P^{1}_{2} - a  & P^{2}_{2} - a^{2}  & P^{3}_{2} - a^{3}   &\cdots& P^{k-1}_{2} - a P^{k-2}_{2}& 1 \\
  P^{1}_{3} - a  & P^{2}_{3} - a^{2}  & P^{3}_{3} - a^{3}   &\cdots& P^{k-1}_{3} - a P^{k-2}_{3}& 1 \\
  \vdots & \vdots & \vdots & \cdots & \vdots &\vdots \\
  P^{1}_{k} - a  & P^{2}_{k} - a^{2}  & P^{3}_{k} - a^{3}   &\cdots& P^{k-1}_{k} - a P^{k-2}_{k}& 1 \\
 \end{vmatrix}
}\\
\end{equation*}

And so on to get,
\begin{equation*}
c \times
\Large{
 \begin{vmatrix}
  P^{1}_{1} - a  & P^{2}_{1} - a^{2}  & P^{3}_{1} - a^{3}   &\cdots& P^{k-1}_{1} - a^{k-1}& 1 \\
  P^{1}_{2} - a  & P^{2}_{2} - a^{2}  & P^{3}_{2} - a^{3}   &\cdots& P^{k-1}_{2} - a^{k-1}& 1 \\
  P^{1}_{3} - a  & P^{2}_{3} - a^{2}  & P^{3}_{3} - a^{3}   &\cdots& P^{k-1}_{3} - a^{k-1}& 1 \\
  \vdots & \vdots & \vdots & \cdots & \vdots &\vdots \\
  P^{1}_{k} - a  & P^{2}_{k} - a^{2}  & P^{3}_{k} - a^{3}   &\cdots& P^{k-1}_{k} - a^{k-1}& 1 \\
 \end{vmatrix}
}\\
\end{equation*}

Multiply $k^{th}$ column with $a^{j}$ and to column j to get,
\begin{equation*}
\norm{X}= c \times
\Large{
 \begin{vmatrix}
  P^{1}_{1}  & P^{2}_{1}   & P^{3}_{1}   &\cdots& P^{k-1}_{1} & 1 \\
  P^{1}_{2}  & P^{2}_{2}   & P^{3}_{2}    &\cdots& P^{k-1}_{2} & 1 \\
  P^{1}_{3}  & P^{2}_{3}   & P^{3}_{3}    &\cdots& P^{k-1}_{3} & 1 \\
  \vdots & \vdots & \vdots & \cdots & \vdots &\vdots \\
  P^{1}_{k}  & P^{2}_{k}   & P^{3}_{k}    &\cdots& P^{k-1}_{k} & 1 \\
 \end{vmatrix}
}\\ 
\end{equation*}

Thus, $ \norm{X} =  (-1)^{k-1} \times c \times \norm{V(P)} \neq 0$

\end{document}